\documentclass{article}

\usepackage{amssymb}

\newbox\smilebox
\newbox\anchorbox
\newbox\noanchorbox
\newbox\tempbox

\setbox\smilebox=\hbox{$\smile$}

\def\anchor{\hbox{\vtop{
           \hbox to \wd\smilebox{\hfil\vrule width.4pt height7pt depth1pt\hfil}
           \vskip  -11.5truept
           \hbox to \wd\smilebox{\hfil$\smile$\hfil}}}}
\setbox\anchorbox=\anchor
\def\noanchor{\hbox{\vtop{
           \hbox to \wd\anchorbox{\hfil\anchor\hfil}
           \vskip -14truept
           \hbox to \wd\anchorbox{\hfil/\hfil}}}}
\setbox\noanchorbox=\noanchor

\def\fg#1#2#3{\setbox\tempbox=\hbox{$\scriptstyle{#2}$}
\ifnum\wd\anchorbox>\wd\tempbox\dimen255=\wd\anchorbox
\else\dimen255=\wd\tempbox\fi
{#1\,\vtop{\hbox to \dimen255{\hfil\anchor\hfil}
           \vskip -6truept
           \hbox to \dimen255{\hfil$\scriptstyle{#2}$\hfil}}
           \,#3}}

\def\nfg#1#2#3{\setbox\tempbox=\hbox{$\scriptstyle{#2}$}
\ifnum\wd\noanchorbox>\wd\tempbox\dimen255=\wd\noanchorbox
\else\dimen255=\wd\tempbox\fi
{#1\,\vtop{\hbox to \dimen255{\hfil\noanchor\hfil}
           \vskip -6truept
           \hbox to \dimen255{\hfil$\scriptstyle{#2}$\hfil}}
           \,#3}}

\setbox1=\hbox{$\bot$}

\def\north#1#2{#1\,
\hbox{$\bot$\llap {\hbox to\wd1 {\hfil $/$\hfil}}}
\,#2}

\def\nao#1#2#3{#1\  \hbox{\vtop{ 
\baselineskip=4pt
\hbox{$\bot$\llap {\hbox to\wd1 {\hfil $/$\hfil}}
\hskip .05em \llap{\hbox{$^{\scriptscriptstyle{a}}$}}}\hbox{$\scriptstyle
{#2}$}}}\, #3}

\def\Ind{\setbox0=\hbox{$x$}\kern\wd0\hbox to 0pt{\hss$\mid$\hss}
\lower.9\ht0\hbox to 0pt{\hss$\smile$\hss}\kern\wd0}

\def\Notind{\setbox0=\hbox{$x$}\kern\wd0\hbox to 0pt{\mathchardef
\nn=12854\hss$\nn$\kern1.4\wd0\hss}\hbox to
0pt{\hss$\mid$\hss}\lower.9\ht0 \hbox to
0pt{\hss$\smile$\hss}\kern\wd0}

\def\ind{\mathop{\mathpalette\Ind{}}}

\def\nind{\mathop{\mathpalette\Notind{}}}

\def\bp{\par{\bf Proof.}$\ \ $}

\def\hbar{\overline{h}}

\def\acl{{\rm acl}}
\def\tp{{\rm tp}}
\def\stp{{\rm stp}}

\def\C{{\frak  C}}

\def\FF{{\bf F}}
\def\K{{\bf K}}
\def\Mod{{\bf Mod}}
\def\Mor{{\bf Mor}}
\def\N{{\mathbb N}}

\def\Q{{\mathbb Q}}
\def\Z{{\mathbb Z}}

\def\tp{{\rm tp}}
\def\stp{{\rm stp}}

\def\dom{{\rm dom}}
\def\acl{{\rm acl}}
\def\dcl{{\rm dcl}}
\def\range{{\rm range}}

\def\Fa0{{\FF^a_{\aleph_0}}}

\def\bp{{\bf Proof.}\quad}
\def\endproof{\medskip}
\def\<{\langle}
\def\>{\rangle}

\newtheorem{Theorem}{Theorem}[section]
\newtheorem{Proposition}[Theorem]{Proposition}
\newtheorem{Definition}[Theorem]{Definition}
\newtheorem{Remark}[Theorem]{Remark}
\newtheorem{Example}[Theorem]{Example}
\newtheorem{Lemma}[Theorem]{Lemma}
\newtheorem{Corollary}[Theorem]{Corollary}

\newtheorem{Question}[Theorem]{Question}

\begin{document}

\title{The Schr\"oder-Bernstein property for $a$-saturated models}


\author{John Goodrick\thanks{Partially supported by a grant from the Facultad de Ciencias, Universidad de los Andes, and by a travel grant from MSRI to attend the BIRS workshop on Neostability Theory (January 29-February 4, 2012).}
\\Departmento de Matem\'aticas\\Universidad de los Andes\\jr.goodrick427@uniandes.edu.co
 \and Michael C.\ Laskowski\thanks{Partially supported
by NSF grant DMS-0901336.}\\Department of Mathematics\\University of Maryland}

\maketitle




\begin{abstract}

A first-order theory $T$ has the {\em Schr\"oder-Bernstein
(SB) property} if any pair of elementarily bi-embeddable models are
isomorphic.  We prove that $T$ has an expansion by constants with
the SB property if and only if $T$ is superstable and non-multidimensional.
We also prove that among superstable theories $T$, the class of a-saturated
models of $T$ has the SB property if and only if $T$ has no nomadic types (see Definition~\ref{nomadic} below).

\end{abstract}

\section{Introduction}

The classical Schr\"oder-Bernstein theorem asserts that if $A$ and $B$ are bi-embeddable sets, i.e., there exist injections $f:A\rightarrow B$
and $g:B\rightarrow A$, then there is a bijection between them.  It is natural to extend this concept to classes $(\K,\Mor)$, where
$\K$ is a class of algebraic structures and $\Mor$ is a distinguished class of injections between elements of $\K$.
We say that $(\K,\Mor)$ has the {\em Schr\"oder-Bernstein (SB) property\/} if any pair of bi-embeddable structures in $\K$ (with respect to
$\Mor$) are isomorphic.  
In this paper, we discuss the SB property for classes $\K$ that are subclasses of $\Mod(T)$, the class of models
of a complete, first-order theory $T$.  Throughout this paper, $\Mor$ will always be taken to be the class of elementary embeddings
(which are necessarily injective).
Thus, we say that a theory $T$ has the SB property if  any two elementarily bi-embeddable models of $T$ are isomorphic.
As examples, the theory of algebraically closed fields of any characteristic has the SB property, but the theory of
dense linear order does not.

%
%
%
One motivation for considering the SB~property is that it should be a nice litmus test for our understanding of the models of $T$: Once we have a sufficiently good understanding of these models (knowing that they are classified by some reasonable collection of invariants, or, conversely, knowing that they are ``wild'' in a suitably precise sense), then we ought to be able to say whether or not $T$ has SB. For instance, by the results of Morley in \cite{Morley67}, if $T$ is countable and $\aleph_1$-categorical (e.g., algbraically closed fields of any characteristic), then the models of $T$ are classified by a single cardinal number invariant (a dimension) which is preserved by elementary embeddings.  This implies that such a $T$ has the SB property. In general, it seems that SB is a fairly strong tameness property, but it is strictly weaker than uncountable categoricity.

Given that we are interested in relating SB to the classification of models of $T$, it is not surprising that we use tools from the so-called \emph{classification theory} of Shelah (now more commonly called \emph{stability theory}), developed in the 1970's and expounded in \cite{bible}. One of the main ideas there was the use of dividing lines amongst theories $T$ (such as superstability, NDOP, NOTOP) to separate those $T$ whose classes of models do admit some kind of classification from those for which this is hopeless. Another idea from \cite{bible} that is very useful for the present paper is the development of a local dimension theory for certain classes of elements within a model using the independence notion known as nonforking.

The Schr\"{o}der-Bernstein property for first-order theories seems to have first been considered by Nurmagambetov in \cite{nur2} and \cite{nur1}, where he showed that for totally transcendental theories it is equivalent to nonmultidimensionality. Various other results around SB were proved in the first author's thesis \cite{my_thesis}, such as:

\begin{Theorem}
\label{nonsuperstable_SB}
If $T$ is not superstable, 
then $T$ does not have SB. Furthermore, if $T$ is unstable, 
then for any cardinal $\kappa$, there is an infinite collection 
of $\kappa$-saturated models of $T$ which are pairwise 
bi-embeddable but pairwise nonisomorphic.
\end{Theorem}

Our previous paper \cite{GL} gives a characterization of which countable weakly minimal theories have the SB property, and this characterization is precise enough to show that for any fixed $T$ the SB property is absolute under forcing extensions of the set-theoretic universe (which does not seem obvious \emph{a priori}). We still do not have a satisfactory characterization of which theories have SB in general.

In the current paper, we address two questions: When is it the case that $T$ has SB after naming a set of constants (which we call ``eventual SB''), and when do the sufficiently saturated models of $T$ have the SB property?

We give a simple complete characterization of eventual SB in Theorem~\ref{eventualtheorem}:
It is equivalent to superstability plus nonmultidimensionality. As for SB for $\kappa$-saturated models, 
with Theorem~\ref{amodelsSB}
we  succeed in characterizing SB for $a$-saturated models when $T$ is superstable.
This result, coupled with
 Theorem~\ref{nonsuperstable_SB}, imply that the only remaining case to consider
is when $T$ is strictly stable.  To extend our methods to such  theories would require
additional knowledge about the prevalence of regular types.

A motivating example to keep in mind is the complete theory $T$ of the additive group $(\Z; +)$. It turns out that $T$ is the theory of all torsion-free abelian groups $G$ such that $[G : pG] = p$ for every prime $p$ (this follows from a more general result by Szmielew in \cite{szmielew}). Hence any model $G$ of $T$ can be decomposed as $G = H \oplus \Q^{\kappa}$ where $H \leq \widehat{\Z}$.  The theory  $T$ is ``classifiable'' according to Shelah's dichotomies.  In fact, $T$  is superstable, non-multidimensional, and weakly minimal; see \cite{prest}). The $a$-saturated 
models of $T$ are simply the models of the form $G = \widehat{\Z} \oplus \Q^{\kappa}$ where the cardinal $\kappa$ is infinite, and it is not hard to see that the class of these models has the SB~property. Furthermore, $T$ has the eventual SB~property since we can add constants for every element in a copy of $\widehat{\Z}$ in some model. However, the class of \emph{all} models of $T$ is more complex, and this does not have the SB~property (this follows from the main theorem of \cite{GL}).

Section~2 introduces some concepts and tools necessary for the main results (countable local pre-weight and low substructures). Along the way, we give a new characterization of the $a$-prime models in any superstable theory (Theorem~\ref{charprime}). Section~3 gives the characterization of SB for $a$-saturated models in a superstable theory, and Section~4 gives the characterization of eventual SB (whose proof depends heavily on our analysis of $a$-saturated models in previous sections).

\textbf{Throughout the paper, we will assume that $T$ is a complete superstable theory unless otherwise noted,} though in a few results we note that $T$ is superstable for emphasis.
Given a complete theory $T$, we always work within
a sufficiently saturated model $\C$ of $T$ (for our purposes, $(2^{|T|})^+$-saturated is enough).  
Our notation is mostly standard and follows \cite{pillay} and \cite{bible}, where the interested reader can find definitions for the terms we use from stability theory (``superstable,'' ``regular type,'' etc.).  As in \cite{pillay}, we call a structure {\em a-saturated} in
place of Shelah's `${\bf F}^a_{\aleph_0}$-saturated,' and we call a model
{\em a-prime\/} instead of `${\bf F}^a_{\aleph_0}$-prime over $\emptyset$.'
When describing dimensions of regular strong types inside models, the notation
$\dim(p,A,M)=\kappa$ means that $p$ is based on $A$ and that
any  maximal $A$-independent set of realizations of $p|A$ inside $M$ has size
$\kappa$.

It is noteworthy that none of the results in this paper have any dependence on the cardinality of the language.

\section{a-prime models of superstable theories}

In this section we focus on the $a$-prime and $a$-saturated models of a superstable theory. Recall (from \cite{pillay}) that a model is $a$-saturated 
if it realizes every strong type over a finite subset, and that a model is $a$-prime just in case it is 
$a$-saturated and embeds into any other $a$-saturated model of its complete theory. 
We will freely use well-known facts about $a$-prime and $a$-saturated models from \cite{pillay} and \cite{bible}.

We first focus on proving a characterization of the $a$-prime models in any superstable theory (Proposition~\ref{charprime}) which is important for our subsequent results. Note that the main use of superstability in the proof of this proposition is the ``ubiquity of regular types.'' In the remainder of the section, we introduce low substructures and some lemmas on dimensions in $a$-saturated models that will be useful later.

Our first definition is a variation on the classical notion of \emph{pre-weight} (see \cite{pillay} or \cite{bible}) which measures the size of a set by how many distinct independent elements can fork with it.

\begin{Definition} \label{countablepreweight}
{\em A set $B$ has {\em countable local pre-weight\/}
if for every  finite  set $A$, every stationary, regular type $p\in S(A)$,
and every $A$-independent set $I\subseteq p(\C)$, there is a countable 
$I_0\subseteq I$ such that $(I\setminus I_0) \ind\limits_A B$.
}
\end{Definition}

\begin{Remark} {\em  
In the definition above,
``local'' refers to the fact that we require that the elements of $I$ come from a single regular type. Classically, weight differs from pre-weight in that the weight is the supremum of the pre-weights of all nonforking extensions. 
Whereas Theorem~\ref{charprime} shows that the a-prime model
of a superstable theory has countable local pre-weight, 
Example~\ref{manyparams} shows that it need not have countable local weight.
}
\end{Remark}

\begin{Lemma}  \label{111}
Suppose that $M$ is any a-prime model and $p\in S(\emptyset)$
is stationary and regular.  For any countable set $A$ and for any
$A$-independent set $I$ of realizations of $p|A$, there is no 
uncountable, pairwise disjoint family $\{E_i:i\in\omega_1\}$ of subsets of
$I$ such that $ E_i \nind\limits_A M$ for each $i$.
\end{Lemma}

\bp  By way of contradiction, suppose an uncountable family 
$\{E_i:i\in\omega_1\}$ existed.  We can clearly assume that each
$E_i$ is finite.  For each $i$ choose a finite tuple $a_i$ from $A$ and
a finite $b_i$ from $M$ such that $ E_i \nind_{a_i} b_i$ for each
$i$.  As $A$ is countable, there is a specific $a^*$ such that $a_i=a^*$
for uncountably many $i$.  Thus, by reindexing, we may assume that
every $a_i=a^*$.  Let $N=M[a^*]$ be the a-prime model over
$M\cup\{a^*\}$.  As $a^*$ is finite, $N$ is also a-prime over $\emptyset$.
Choose a maximal, independent set $J\subseteq p(N)$.  As $N$ is 
a-prime, $J$ is countable.  

Next, for each $i$, choose a finite $J(i)\subseteq J$ so that
$a^*b_i \ind_{J(i)} J$.  Arguing as above, we may assume that
there is a single $J^*$ so that $J(i)=J^*$ for every $i$.
Furthermore, since $J^*$ is finite and the sets
$\{E_i\}$ are independent, by eliminating at most finitely many $i$
we may additionally assume that $E_i \ind J^*$
for each $i$.

Now we obtain a contradiction by fixing any remaining $i$.
As $N$ is a-saturated, we can choose $E_i'\subseteq N$
such that $\stp(E_i'a^*b_iJ^*)=\stp(E_ia^*b_iJ^*)$.
Since every element of $E_i'$ realizes $p(N)$ and is independent of
$J^*$, the maximality of $J$ implies that $J$ dominates $E_i'$ over $J^*$.
Since $E_i$ is independent from $J^*$ and forks with $a^*b_i$ over
$\emptyset$, it follows that $E_i$, and hence $E_i'$, forks with $a^*b_i$
over $J^*$.  Combining this with the  domination described above implies that
$ a^*b_i \nind_{J^*} J$, which is a contradiction.
\endproof

\begin{Lemma}  \label{222}
Suppose that $M$ is any a-prime model and $p\in S(\emptyset)$
is stationary and regular.  For any countable set $A$ and for any
$A$-independent set $I$ of realizations of $p|A$, there is a countable
$I^*\subseteq I$ such that
 $ (I\setminus I^*) \ind\limits_{AI^*} M$.
\end{Lemma}

\bp  We first argue that for any countable set $A$ there is a countable
$I_0\subseteq I$ such that $(I\setminus I_0) \ind_A M$.
To see this, given a countable set $A$, call a subset $E\subseteq I$
a {\em minimal witness to forking} if $E \nind_A E M$, but any proper
subset of $E$ is free from $M$ over $A$.  It is clear that every minimal
witness is finite, and that if we set $I_0$ to be the union of all the minimal
witnesses, then $(I\setminus I_0) \ind_A M$.  Thus, it remains to prove
that $I_0$ is countable.

However, if $I_0$ were uncountable, then we would have uncountably
many minimal witnesses $\{E_i\}$.  By the $\Delta$-system lemma,
there would be a finite set $G\subseteq I$ and an uncountable family
$\{E_i:i\in\omega_1\}$ such that $E_i\cap E_j=G$ for distinct $i,j$.
But then, apply Lemma~\ref{111} with $A'=A\cup G$, $I'=I\setminus G$,
and the family $\{F_i\}$, where $F_i=E_i\setminus G$ and obtain a 
contradiction.

Now, to prove the Lemma, suppose we are given
a countable set $A$.  Form a sequence $I_0\subseteq I_1\subseteq \dots$
of  countable subsets of $I$ by applying the result in the first paragraph
successively to the  countable sets $A$, $A\cup I_0$,  $A\cup  I_1$,
et cetera.  Then the set $I^*=\bigcup I_n$ satisfies our demands.
\endproof

\begin{Proposition} \label{preweight}
 Every a-prime model has countable local pre-weight.
In fact, given any sets $B\subseteq A$ with $B$ finite and $A$
countable, and for any stationary, regular $p\in S(B)$, then for every 
$B$-independent set $I\subseteq p(\C)$, there is a countable set
$I^*\subseteq I$ such that $(I\setminus I^*) \ind\limits_{AI^*} M$.
\end{Proposition}

\bp  First, as $A$ is countable, there is a finite $I_0\subseteq I$
such that $I\setminus I_0$ is $AI_0$-independent.  
Thus, by replacing $I$ by $I\setminus I_0$ and $A$ by $A\cup I_0$,
we may additionally assume that $I$ is $A$-independent.

Next, let $N=M[B]$ be a-prime over $M\cup B$.  Then $N$ is also a-prime.
But furthermore, as $B$ is finite, $N$ is also a-prime over $B$.
Thus, if we work over $B$ and apply Lemma~\ref{222}, we obtain
the requisite $I^*$.
\endproof

In fact, the countability of local pre-weight characterizes the a-prime models among
the class of all a-saturated models.

\begin{Proposition} \label{charprime}
 The following are equivalent for an a-saturated model $M$ of a 
superstable theory:
\begin{enumerate}
\item  $M$ is a-prime over $\emptyset$;
\item  Every infinite, indiscernible set in $M$ is countable;
\item  For all finite $B\subseteq M$, every stationary, regular $p\in S(B)$
has countable dimension;
\item  $M$ has countable local pre-weight.
\end{enumerate}
\end{Proposition}

\bp  The equivalence of (1) and (2) is given in IV~4.18 of \cite{bible},
noting that any finite tuple is trivially ${\bf F}^a_{\aleph_0}$-atomic
over $\emptyset$.   (2) implies (3) is
immediate.  

To see that (3) implies (2), suppose that there is an uncountable
indiscernible set $I\subseteq M$.  Let $q=Av(I,M)$.  By superstability there
is a regular type $p$ non-orthogonal to $q$.  Since $M$ is a-saturated, possibly
by replacing $p$ by a non-orthogonal regular type, we may assume
that $p$ is based and stationary on a finite $B\subseteq M$.  Again by
superstability, there is a finite $I_0\subseteq I$ on which $q$ is stationary
and moreover, by padding $B$ with a finite Morley sequence in $p$,
we may additionally assume that the types $p',q'\in S(BI_0)$, that are
parallel to $p,q$ respectively, are not almost orthogonal.  From this and the
fact that $p'$ has weight one, it is clear that $M$ must contain an uncountable
Morley sequence in $p'$.

Finally, (4) implies (3) is obvious, and the implication (3) implies (4) is the
content of Proposition~\ref{preweight}.
\endproof

\begin{Definition} {\em  Given an a-saturated
model $N$, $M$ is a {\em low substructure
of $N$\/} if  $M\preceq N$, 
$M$ is a-prime, and $\dim(p,M,N)\ge\aleph_0$ for every
regular  type $p\in S(M)$.
}
\end{Definition}

\begin{Lemma}  \label{existence}
Every a-saturated model has a low substructure.
\end{Lemma}

\bp  It suffices to prove that every a-prime model
has a low substructure.
By the uniqueness of a-prime models, it suffices to construct
a single  a-prime model  $N$ that has a low
substructure.    

Toward this end, fix $M$ any a-prime model.
Let $\Gamma\subseteq S(M)$ be any maximal subset of pairwise
orthogonal weight one types over $S(M)$.
Let $$I=\bigcup_{p\in\Gamma} I_p$$
be independent over $M$ such that each $I_p$ is a Morley sequence
of length $\omega$ built from $p$, and let $N$ be a-prime
over $M\cup I$.  By construction, $\dim(q,M,N)\ge\aleph_0$ for
every $q\in S(M)$, so it
suffices to show that $N$ is a-prime over $\emptyset$.
By Lemma~\ref{charprime}, it suffices to show that as a set,
$N$ has countable local pre-weight.
To see this, choose any finite set $A$, any stationary, regular $r\in S(A)$,
and any Morley sequence $J$ built from $r$.  There is at most one
$p\in\Gamma$ non-orthogonal to $r$.  Choose $M'\preceq N$
to be a-prime over $M\cup I_p$ if such a $p\in\Gamma$ exists,
or else let $M'=M$ if $r\perp p$ for every $p\in\Gamma$.  
Note that in either case, $\tp(N/M')$ is orthogonal to $r$.
Since $M$ is a-prime it has countable local pre-weight,
so there is a countable $J_0\subseteq J$
with $M \ind_{AJ_0} J$.  By the construction of $M'$, it follows
that there is a countable $J_1$, $J_0\subseteq J_1\subseteq J$, satisfying
 $M' \ind_{AJ_1} J$.  As $\tp(N/M')$ is orthogonal to $r$,
$ N \ind_{AJ_1} J$, so $N$ has countable local pre-weight.
\endproof

\begin{Lemma}  \label{Xmatch}
Suppose that $N$ is an a-saturated
model, $X\subseteq N$
is any set, and $p,q\in S(X)$ are  stationary, regular
types that are not almost orthogonal.  Then $\dim(p,X,N)=\dim(q,X,N)$.
\end{Lemma}

\bp  
We first prove that for each $c\in p(N)$, there is $d\in q(N)$ that forks
with $c$ over $X$.  To see this, by non-almost orthogonality, for any
such $c$ there is $d_0\in q(\C)$  forking with $c$ over $X$. 
Choose $B\subseteq X$ finite such that $p$ and $q$ are based on $B$
and $c \nind_B d_0$.  As $N$ is a-saturated, we can find $d\in N$
such that $\stp(d/Bc)=\stp(d_0/Bc)$.  To see that $d$ realizes $q$,
it suffices to show that $d \ind_B X$.  However, if it forked, then $r=\stp(d/X)$
would be a forking extension of a strong type parallel to $q$.  But then,
by regularity $r$ would be orthogonal to $p$, which would be a contradiction.

Now let
 $I\subseteq p(N)$ be any maximal $X$-independent set.
From the argument above, for each $c\in I$, choose $d_c\in q(N)$
that forks with $c$ over $X$.  It follows immediately from the fact
that regular types have weight one that the mapping $c\mapsto d_c$ is
injective, and moreover that $J=\{d_c:c\in I\}$ is $X$-independent.
Thus, $\dim(p,X,N)\le\dim(q,X,N)$.  By symmetry, this suffices to prove
the lemma.
\endproof

\begin{Proposition} \label{dimensionmatch}
 Suppose $M_0,M_1$ are both low substructures
of an arbitrary a-saturated
model $N$.  If $p\in S(M_0)$ and $q\in S(M_1)$ are 
non-orthogonal regular types, then $\dim(p,M_0,N)=\dim(q,M_1,N)$.
\end{Proposition}

\bp  By the definition of a low substructure,
$M_0$ and $M_1$ are both a-prime and each dimension is infinite.  If the dimensions
are both countable, they are equal.  Thus, by symmetry, assume that
$\dim(p,M_0,N)=\kappa>\aleph_0$.  It suffices to show that
$\dim(q,M_1,N)=\kappa$.  To see this, first choose a finite $B\subseteq M_0$
on which $p$ is based and stationary.  As $M_0$ is a-prime, $\dim(p,B,M_0)$
is countable, hence $\dim(p,B,N)=\kappa$ as well.  Let $I\subseteq N$
be a maximal $B$-independent set of realizations of $p|B$.  
By Proposition~\ref{preweight} there is a countable $I_0\subseteq I$
such that $ (I\setminus I_0) \ind_{BI_0} {M_1}$.  By adding at most finitely
many additional points to $BI_0$, we may assume that the types parallel
to $p$ and $q$ over $M_1BI_0$ are not almost orthogonal.  Thus, by 
Lemma~\ref{Xmatch},
as $|I\setminus I_0|=\kappa$, it follows
that $\dim(q,M_1BI_0,N)=\kappa$.  But then, as $BI_0$ is countable,
it follows that $\dim(q,M_1,N)=\kappa$ and we finish.
\endproof

\begin{Lemma}  \label{preserved} $T$ superstable.
If $f:M\rightarrow N$ is any elementary embedding of a-saturated models,
and if $M_0$ is a low substructure of $M$, then $f(M_0)$ is a low
substructure of $N$.
\end{Lemma}

\bp  Clearly, $f(M_0)\preceq N$ and is a-prime, since it is isomorphic to $M_0$.
Furthermore, if $q\in S(f(M_0))$ is regular, then as $p:=f^{-1}(q)$
is a regular type over $M_0$, $p$ has infinite dimension in $M$, hence
$q$ has infinite dimension in $N$.  Thus, $f(M_0)$ is a low substructure
of $N$.
\endproof

\begin{Corollary}  \label{wrap} 
Suppose that $T$ is superstable,
$M_0$ is a low substructure of an a-saturated model $M$,
and that $f:M\rightarrow M$ is an elementary endomorphism.
If $p\in S(M_0)$ is regular and non-orthogonal to $f(p)$, 
then $\dim(p,M_0,M)=\dim(f(p),f(M_0),M)$.
\end{Corollary}

\bp This follows immediately from Proposition~\ref{dimensionmatch}
and Lemma~\ref{preserved}.
\endproof

We close this section with a remark about Definition~\ref{countablepreweight}.
For the moment, say that a set $B$ has {\em countable local weight}
if for any set $A$ independent from $B$, 
any stationary, regular $p\in S(A)$, and any
$I\subseteq p(\C)$, there is a countable $I_0\subseteq I$ such that
$ (I\setminus I_0) \ind_A B$.  The following example shows that having countable
local weight is
too much to expect for the universe of an a-prime model, even if the
theory is weakly minimal and unidimensional.  

\begin{Example}  \label{manyparams} The a-prime model of
$Th((\Z,+))$ does not have countable local weight.  In fact, the a-prime model
$(G,+,\dots)$ of any weakly minimal
group  in a countable language that has a family of $2^{\aleph_0}$ strong
types over $\emptyset$, each describing cosets of the principal generic
type $p$ does not have countable local weight.
\end{Example}

\bp Let $(G,+,\dots)$ be the a-prime model of such a theory.  As the language
is countable, we can inductively find a sequence 
$\{q_\alpha:\alpha\in\omega_1\}$ of strong types that are `almost independent
over $\emptyset$' i.e., for any choices of $b_\alpha$ realizing $q_\alpha$,
the set $B=\{b_\alpha:\alpha\in\omega_1\}$ is independent over $\emptyset$.

As $G$ is a-saturated, we can choose such a set $B\subseteq G$
as described above.
Let $A=\{a_\alpha:\alpha\in\omega_1\}$ be
free from $G$ over $\emptyset$, with each $a_\alpha$ realizing $q_\alpha$.
Now let $I=\{c_\alpha:\alpha\in\omega_1\}$, where each 
$c_\alpha:=b_\alpha-a_\alpha$ realizes the principal generic type $p$.
It is easy to see that $I$ is $A$-independent, but $I_1 \nind_A G$ for any
non-empty $I_1\subseteq I$.
\endproof

\section{The SB-property for a-saturated models}

In this section, we use the results of the previous section to characterize which superstable theories have the SB property for $a$-saturated models: they are precisely the ones without nomadic types (see Definition~\ref{nomadic} and Theorem~\ref{amodelsSB}). It turns out that this condition is slightly stronger than being non-multidimensional (nmd), and so we found it useful to first establish some facts about $a$-prime models in any superstable nmd theory.

We continue to assume that $T$ is superstable unless otherwise specified.

\begin{Definition}  \label{nomadic}
{\em
A non-algebraic strong type $p$ is {\em nomadic\/} if 
there is an automorphism $f\in Aut(\C)$
such that the $n$-fold iterate $f^{(n)}(p)\perp p$ for each integer $n$.
}
\end{Definition}

It is easy to see that if a superstable theory $T$ has a nomadic type, then
it has a regular nomadic type.  As well, it is easily seen (using e.g.,
Lemma~1.4.3.3 of \cite{pillay}) that any type $p$ that is orthogonal
to $\emptyset$ is nomadic, so that any multidimensional theory
has nomadic types.  However, even some nmd theories have nomadic types: 

\begin{Example}
Let $T$ be the theory of a collection $\{E_i : i \in \N\}$ of refining equivalence relations ($E_{i+1}(x,y) \Rightarrow E_i(x,y)$) such that $E_0$ has exactly two classes and each $E_i$-class splits into two $E_{i+1}$-classes. Then $T$ has quantifier elimination and is complete, superstable (in fact, weakly minimal), and nmd, and $T$ has a unique $1$-type $p(x)$ over $\emptyset$. If $f \in \textup{Aut}(\C)$ is chosen so that it permutes the $E_i$-classes in a single cycle of order $2^{i+1}$, then $p \perp f(p) \perp f^2(p) \perp \ldots$, so $p$ is nomadic.
\end{Example}

Among nmd theories,
then the existence of nomadic types can be obviated by eliminating
automorphisms of $\acl(\emptyset)$:

\begin{Lemma}  Suppose $T$ is nmd.  If,
in $\C^{\rm eq}$, $\acl(\emptyset)=\dcl(\emptyset)$, then
$T$ has no nomadic types.
\end{Lemma}

\bp    By the observation
above, it suffices to show that no stationary, regular $p\in S(A)$
is nomadic via any
automorphism $f\in Aut(\C)$.
In fact, we show that $p\not\perp f(p)$ for any such $p$ and $f$.
To see this, let $A'=f(A)$, and choose $A''$ having the same strong
type as $A$ over $\emptyset$, but with $A''$ free from $A\cup A'$ over
$\emptyset$.  Note that by our assumptions,
$\tp(A/\acl(\emptyset))=\tp(A''/\acl(\emptyset))=\tp(A'/\acl(\emptyset))$.
Let $p''$ be the (regular) type conjugate to $p$ over
$A''$.  Because $T$ has nmd, $p\not\perp \emptyset$, so by
e.g., Lemma~1.4.3.3 of \cite{pillay}, $p\not\perp p''$ and
$f(p)\not\perp p''$.  Thus, $p\not\perp f(p)$ by the transitivity of
nonorthogonality among regular types.
\endproof

\begin{Lemma}  \label{nmd3model}  Suppose that $T$ is
nmd and $M\preceq M^*\preceq N$, where $M$ and $N$ are
a-saturated and $M^*$ is any model that is not equal to $N$. Then there is $c\in N\setminus M$ such that $\tp(c/M)$ is regular
and $c \ind\limits_M N$.
\end{Lemma}

\bp 
By superstability and the fact that
$M^*\neq N$, there is a regular $q\in S(M^*)$ 
that is realized in $N$ (see e.g., Proposition~8.3.2 of \cite{pillay}).
Since $T$ has nmd,
$q$ is not orthogonal to $M$.  So, since $M$ is a-saturated,
the existence of such a $c$ follows immediately by Proposition~8.3.6 of
\cite{pillay}.
\endproof

\begin{Proposition}  \label{primeoveramodels}
Suppose that $T$ is nmd and that $N \models T$ is any model with an a-saturated elementary submodel $M$.
Then $N$ is a-saturated.
\end{Proposition}

\bp  Let $N^*=M[N]$ be any a-prime model over $M\cup N=N$.  If
$N\neq N^*$, then by Proposition~\ref{nmd3model} there would be 
$c\in N^*\setminus M$ with $c \ind_M N$, which contradicts the fact
that $N^*$ is dominated by $N$ over $M$.
\endproof

Thus, for any set $B$ that contains an a-saturated model $M$,
any model $N$ containing $B$ is automatically a-saturated.
Hence,
the notions of `prime over $B$' and `minimal over $B$' are equivalent
to the notions `a-prime over $B$' and `a-minimal over $B$', respectively.  
These observations will be used extensively in the next 
section.\footnote{These equivalences also illustrate an error in popular
parlance.  In the setting of superstable, NDOP theories in a countable language,
in Chapter~12 of \cite{bible} Shelah proved that NOTOP is equivalent
to the `$(\infty,2)$-existence property'.  In several places, this
unweildly phrase has been replaced by `PMOP', an acronym for `Prime Models
Over Pairs.'  Although this term sounds better, it is misleading.
The results above show that for any superstable, nmd theory (even those
with OTOP, see e.g., Example~2.2 of \cite{HL}), 
there is a prime model over $M_1\cup M_2$, where $M_1$ and $M_2$
are any models containing an a-saturated model.  
The `correct' meaning of PMOP (i.e., equivalent to $(\infty,2)$-existence)
is that there is a $t$-constructible model
over any independent pair of models.}

\begin{Proposition}  \label{structure}
Suppose that $T$ is nmd and that $M\preceq N$ are both a-saturated and $J$ is any
maximal, $M$-independent
subset of $N$ consisting of realizations of regular types over $M$.
Then $N$ is a-prime and a-minimal over $M\cup J$.  
\end{Proposition}

\bp  
Let $M^*\preceq N$ be any a-prime model over $M\cup J$.
To see that $N$ is both a-prime and a-minimal over $M\cup J$, it
suffices to prove that $M^*=N$.  If this were not the case, then
by Lemma~\ref{nmd3model} there would be $c\in N\setminus M^*$
such that $\tp(c/M)$ is regular and $c \ind_M {M^*}$, which would contradict
the maximality of $J$.
\endproof

\begin{Definition}  \label{dimpreserving}
{\em 
Suppose that $M_0\preceq M$ and $N$ are a-saturated.
An elementary embedding $f:M\rightarrow N$ is {\em dimension preserving
over $M_0\preceq M$} if $\dim(p,M_0,M)=\dim(f(p),f(M_0),N)$ for
every regular $p\in S(M_0)$.  
}
\end{Definition}

\begin{Corollary}  \label{criterion}  
Suppose $T$ is nmd.
If
there is an elementary embedding
$f:M\rightarrow N$ between a-saturated
models that is dimension-preserving over some
a-saturated $M_0\preceq M$, then $M\cong N$.
In fact, the isomorphism can be chosen to extend $f|_{M_0}$.
\end{Corollary}

\bp  Fix a-saturated models $M, N$ and 
an elementary embedding $f:M\rightarrow N$ that is dimension preserving
over some a-saturated  $M_0\preceq M$.  
Let $\Gamma\subseteq S(M_0)$ be a maximal, pairwise orthogonal
set of regular types over $M_0$, and for each $p\in\Gamma$, let $I_p$ be
a maximal, independent subset of $p(M)$.  Note that 
$I=\bigcup_{p\in\Gamma}I_p$ is a maximal, $M_0$-independent set of
realizations of regular types in $M$.  
Next, for each $p\in \Gamma$, let $J_p$ be a maximal, $f(M_0)$-independent
set of realizations of $f(p)$ in $N$.
As $f|_{M_0}$ is an isomorphism between $M_0$ and $f(M_0)$, it follows
that $J=\bigcup_{p\in\Gamma} J_p$ is a maximal, $f(M_0)$-independent 
set of realizations of regular types in $N$.
Since $f$ is dimension preserving over $M_0$, $|I_p|=|J_p|$ for each
$p\in \Gamma$.  For each $p\in\Gamma$, fix any bijection
$g_p:I_p\rightarrow J_p$.   By indiscernibility, for each $p\in\Gamma$,
the map $h_p:=f|_{M_0}\cup g_p$  is elementary, and by independence over
a model, $h:=\bigcup h_p:M_0\cup I\rightarrow f(M_0)\cup J$ is elementary
as well.   By Proposition~\ref{structure}, $M$ is a-prime over 
$\dom(h)$, while $N$ is a-prime over $\range(h)$.  It follows
from the uniqueness of a-prime models that $h$ extends to an
isomorphism $h^*:M\rightarrow N$.
\endproof

The proof of the following Lemma is immediate.

\begin{Lemma}  \label{monotone}  
If $f:M\rightarrow N$ is any elementary embedding and $M_0\preceq M$
is arbitrary, then $\dim(p,M_0,M)\le\dim(f(p),f(M_0),N)$ for any regular
type $p\in S(M_0)$.
\end{Lemma}

\begin{Proposition}  \label{matching}
Assume that $T$ has no nomadic types.
Suppose that $M$ and $N$ are both a-saturated and 
$f:M\rightarrow N$ and $g:N\rightarrow M$ are elementary embeddings.
Then $f$ is dimension preserving over any
low substructure $M_0$ of  $M$.
\end{Proposition}

\bp  Let $h=g\circ f$ denote the composition.  Fix any low substructure
$M_0$ of
$M$ and any regular type $p\in S(M_0)$.  Since there are no nomadic
types, there is a positive integer $n$ so that  the $n$-fold
composition $k=h^{(n)}$ satisfies $p\not\perp k(p)$.
By Corollary~\ref{wrap} we have $\dim(p,M_0,M)=\dim(k(p),k(M_0),M)$.
But, by iterating Lemma~\ref{monotone} we have
$$\dim(p,M_0,M)\le\dim(f(p),f(M_0),N)\le\dim(k(p),k(M_0),M)$$ 
so $\dim(p,M_0,M)=\dim(f(p),f(M_0),N)$
as required.
\endproof

\begin{Theorem}  \label{amodelsSB} For a superstable theory $T$, the following are equivalent:
\begin{enumerate}
\item $T$ has the Schr\"oder-Bernstein property for $a$-saturated models;
\item there is no \emph{infinite} collection of pairwise elementarily bi-embeddable, pairwise nonisomorphic $a$-saturated models of $T$;
\item $T$ has no nomadic types.
\end{enumerate}
\end{Theorem}

\bp  $1 \Rightarrow 2$ is trivial. The direction $2 \Rightarrow 3$  is proved in Theorem~4.8 of \cite{GL} (whose statement does not mention saturation, but as noted in the proof there, the argument can be used to produce bi-embeddable, nonisomorphic $a$-saturated models).
Finally, for $3 \Rightarrow 1$, note that $T$ superstable with no nomadic types
implies that 
$T$ is nmd as well.  Choose a-saturated  models $M$ and $N$ and fix elementary embeddings $f:M\rightarrow N$ and $g:N\rightarrow M$.
By Lemma~\ref{existence}, choose  a low substructure $M_0$ of $M$.
By Proposition~\ref{matching}, $f$ is dimension preserving over $M_0$,
so $M$ and $N$ are isomorphic by Corollary~\ref{criterion}.
\endproof

It might be true that having no nomadic types implies the SB~property for $a$-saturated models in \emph{any} stable theory; we know of no counterexample. (For a strictly stable $T$, ``$a$-saturated'' means ``all strong types over subsets of size less than $\kappa_r(T)$ are realized.'')

\section{The eventual SB property for models}

\begin{Definition} 
{\em   
A complete theory $T$ has the {\em eventual SB property\/} if there is
a \emph{small} set $A\subseteq \C$ such that the expansion
$Th_A(\C)$ formed by adding a new constant symbol for each element of
$A$ has the SB property. (Here ``small'' means that $|A| < \kappa$ for some cardinal $\kappa$ such that the universal domain $\C$ is $\kappa$-saturated.)
}
\end{Definition}

The goal of this short section is to characterize those theories with
the eventual SB property.



\begin{Theorem}  \label{eventualtheorem}
The following are equivalent for a complete theory $T$:
\begin{enumerate}
\item  $T$ has the eventual SB property;
\item  For every small subset $A \subseteq \C$ containing an a-saturated model, 
$Th_A(\C)$ has the SB property;
\item  $T$ is superstable and nmd.
\end{enumerate}
\end{Theorem}

\bp  $2\Rightarrow 1$:  Trivial.

$1\Rightarrow 3$:  We prove the contrapositive. If $T$ is \emph{not} superstable, then the same is true of $Th_A(\C)$ for any small set $A \subseteq \C$, so $Th_A(\C)$ does not have SB by Theorem~\ref{nonsuperstable_SB}. The other case to consider is when $T$ is stable and multidimensional, in which case $Th_A(\C)$ is again stable and multidimensional, and the failure of SB follows from Theorem~4.8 of \cite{GL} (noting that any regular type $p$ orthogonal to $\emptyset$ satisfies the hypothesis of that result).

$3\Rightarrow 2$:  Fix a small set $A$ containing an a-saturated model $M$,
and choose any pair $N_1^*, N_2^*$ of bi-embeddable models $Th_A(\C)$. 
That is, the reducts $N_1$ and $N_2$ to the original language are models
of $T$ that are bi-embeddable over $A$.  
By Proposition~\ref{primeoveramodels}, both reducts $N_1$ and $N_2$ are
themselves a-saturated.  We argue that any  $L(A)$-elementary embedding
$f:N_1^*\rightarrow N_2^*$, when viewed as an $L$-elementary embedding
from $N_1$ to $N_2$ that fixes $A$ pointwise, is dimension preserving
over any $M_1\preceq N_1$ that is a-prime over $A$.
To see this, fix any such function $f$ and choose any  $M_1\preceq N_1$
that is a-prime over $A$.   
As $f$ is over $A$,  $f(M_1)$, which we denote
by $M_2$, is also a-prime over $A$.

Note that for each regular type $q\in S(M)$ and each $i=1,2$,
the fact that $M_i$ is dominated by $A$ over $M$ implies that
the non-forking extension $q|A$ is omitted in $M_i$.
Furthermore, as $q|A$ is fixed
by any embedding over $A$, Lemma~\ref{monotone} and the fact that
$N_1^*$ and $N_2^*$ are bi-embeddable $L(A)$-structures imply
that $\dim(q,A,N_1)=\dim(q,A,N_2)$.    It follows that 
$\dim(q,M_1,N_1)=\dim(q,M_2,N_2)$ for every regular $q\in S(M)$.
However, it follows from nmd and the fact that $M$ is an a-saturated
model that
any regular type $p\in S(M_1)$ is non-orthogonal to some regular
$q\in S(M)$.  Moreover, $f(p)$ is non-orthogonal to the same type $q$.
Thus, applying Lemma~\ref{Xmatch} on both sides yields
$$\dim(p,M_1,N_1)=\dim(q,M_1,N_1)=\dim(q,M_2,N_2)=\dim(f(p),M_2,N_2)$$
Thus, $f$ is dimension-preserving over $M_1$.  By 
Corollary~\ref{criterion}, the $L$-structures $N_1$ and $N_2$ are isomorphic
via an isomorphism $h$ extending $f|A$.  As $f$ fixes $A$ pointwise,
$h$ is an $L(A)$-isomorphism between $N_1^*$ and $N^*_2$.
\endproof

\begin{Question}
What conditions on a set $A$ are needed to ensure that $Th_A(\C)$ has SB for a
classifiable, nmd theory?
\end{Question}

\bibliographystyle{amsplain}

\bibliography{modelth}

\end{document}